\documentclass{amsart} 
\usepackage{amssymb}
\usepackage{amsmath}
\usepackage{amsfonts}

\begin{document}
\newtheorem{theo}{Theorem}[section]
\newtheorem{prop}[theo]{Proposition}
\newtheorem{lemma}[theo]{Lemma}
\newtheorem{coro}[theo]{Corollary}
\theoremstyle{definition}
\newtheorem{exam}[theo]{Example}
\newtheorem{defi}[theo]{Definition}
\newtheorem{rem}[theo]{Remark}


\newcommand{\Bb}{{\bf B}}
\newcommand{\Cb}{{\bf C}}
\newcommand{\Nb}{{\bf N}}
\newcommand{\Qb}{{\bf Q}}
\newcommand{\Rb}{{\bf R}}
\newcommand{\Zb}{{\bf Z}}
\newcommand{\Ac}{{\mathcal A}}
\newcommand{\Bc}{{\mathcal B}}
\newcommand{\Cc}{{\mathcal C}}
\newcommand{\Dc}{{\mathcal D}}
\newcommand{\Fc}{{\mathcal F}}
\newcommand{\Ic}{{\mathcal I}}
\newcommand{\Jc}{{\mathcal J}}
\newcommand{\Kc}{{\mathcal K}}
\newcommand{\Lc}{{\mathcal L}}
\newcommand{\Mx}{{\mathcal M}}
\newcommand{\Oc}{{\mathcal O}}
\newcommand{\Pc}{{\mathcal P}}
\newcommand{\Sc}{{\mathcal S}}
\newcommand{\Tc}{{\mathcal T}}
\newcommand{\Uc}{{\mathcal U}}
\newcommand{\Vc}{{\mathcal V}}

\author{Nik Weaver}

\title [Set theory and C*-algebras]
       {Set theory and C*-algebras}

\address {Department of Mathematics\\
          Washington University in Saint Louis\\
          Saint Louis, MO 63130}

\email {nweaver@math.wustl.edu}

\date{\em April 8, 2006}

\begin{abstract}
We survey the use of extra-set-theoretic hypotheses, mainly the
continuum hypothesis, in the C*-algebra literature. The Calkin
algebra emerges as a basic object of interest.
\end{abstract}

\maketitle


Recently Charles Akemann and the author \cite{AW} used Jensen's diamond
principle to solve an old open problem about the existence of nontrivial
C*-algebras with only one irreducible representation. This raises
the question: in general, to what extent is extra-set-theoretic
reasoning relevant to C*-algebras?

The one extra axiom that has been repeatedly used
by C*-algebraists is the continuum hypothesis, although there are
a few instances where it was noted that Martin's axiom or some
more special weakening of CH would suffice to settle the problem
at hand. Interestingly, it appears that C*-algebraists generally
tend to regard a problem as solved when it has been answered using
CH. This may have to do with the fact that in most cases
the other direction of the presumed independence result would involve
set theory at a substantially more sophisticated level. There could
be an opportunity for set theorists here; it seems likely that most
of the theorems we survey below are the easy halves of independence
results, the hard --- or more set-theoretically sophisticated
--- directions of which have not yet been proven.

CH has been particularly valuable in the study of the Calkin algebra
$\Cc(l^2) = \Bc(l^2)/\Kc(l^2)$ (see Example 1.5 below). This algebra can
be viewed as a noncommutative or ``quantum'' analog of the Stone-\v{C}ech
remainder $\beta\omega - \omega$ \cite{vM}. It should be of basic
set-theoretic interest and a rich source of future work.

In the first section below I give a quick introduction to C*-algebras.
In the second section I survey a variety of consistency results taken
from the C*-algebra literature, and in the third I focus on the Calkin
algebra. My concern is consistency results and I do not discuss other
relations between C*-algebras and set theory or logic, such as the
Glimm-Effros dichotomy \cite{HKL} which arose out of C*-algebra
representation theory, connections with nonstandard analysis \cite{AGPS, HO},
Mundici's work on decidability of isomorphism of AF-algebras
\cite{Mun}, or simply the use of basic set-theoretic techniques that
do not go beyond ZFC (e.g., my solution of Dixmier's problem using a
transfinite recursion of length $2^{\aleph_0}$ \cite{W2}).

I wish to thank Ilijas Farah and the referee for providing several
suggestions for improving the exposition of this paper.

\section{C*-algebras}\label{sect1}

We begin with a brief introduction to C*-algebras. Some good general
references are \cite{Con, KR, W1}, or, at a more technical level,
\cite{Ped, Sak, Tak}. Our set-theoretic terminology is standard
\cite{Jec, Kun}.

\subsection{Basic definitions}

Let $H$ be a complex Hilbert space. The scalar field will be complex
throughout.

There is an inconsistency in the mathematics and physics literature
as to whether the inner product $\langle \cdot, \cdot\rangle$ on $H$
should be linear in the first variable and antilinear in the second,
or vice versa. We follow the mathematical convention which takes the
inner product to be linear in the first variable. Thus, $\Cb^n$ is
a Hilbert space with inner product $\langle v,w\rangle =
\sum_1^n v_i\bar{w}_i$ and $l^2$ is a Hilbert space with
inner product $\langle v,w\rangle = \sum_1^\infty v_i\bar{w}_i$.

A {\it bounded operator} on $H$ is a linear map $A: H \to H$ whose
{\it norm}
$$\|A\| = \sup\left\{\frac{\|Av\|}{\|v\|}: v \neq 0\right\}$$
is finite, where $\|v\| = \sqrt{\langle v,v\rangle}$. Here we
adopt the usual notation according to which $A(v)$
is written $Av$. This originates in the finite-dimensional case, where
$H \cong \Cb^n$ and $A$ can be identified with an $n \times n$ matrix
$\tilde A$, with $A(v)$ being the matrix product of $\tilde{A}$ with $v$.

Let $\Bc(H)$ be the set of all bounded operators on $H$. Then $\Bc(H)$
is a Banach space with the obvious linear structure and the norm defined
above. Moreover, it carries two additional algebraic operations, a
product and an involution, which all
together give it a very robust structure. The product of two operators
is defined to be their composition and the involution on $\Bc(H)$ is
the bijection $*: \Bc(H) \to \Bc(H)$ that takes an operator
$A$ to its {\it adjoint} $A^*$, which is characterized by the condition
$$\langle Av, w\rangle = \langle v, A^*w\rangle$$
for all $v,w \in H$. When $H = \Cb^n$ and $A$ is identified with
an $n \times n$ matrix $(a_{ij})$, its adjoint has matrix
$(\bar{a}_{ji})$.

An operator $A$ is {\it self-adjoint} if $A = A^*$. These operators
play the role of ``real'' (as opposed to complex) elements of $\Bc(H)$,
primarily because of the spectral theorem for bounded self-adjoint
operators. There are several versions of this result, but the
easiest to state and probably the most useful is the following: if
$A \in \Bc(H)$ is self-adjoint then there is a measure space
$(X,\mu)$, a bounded measurable real-valued function $f$ on $X$, and a
Hilbert space isomorphism $U: H \cong L^2(X,\mu)$ taking $A$ to
$M_f$, the operator of multiplication by $f$. That is, $M_f$
is the operator on $L^2(X,\mu)$ defined by $M_f(g) = fg$, and
we have $A = U^{-1}M_fU$. Thus a self-adjoint operator on $H$ ``is''
a realization of $H$ as some $L^2(X,\mu)$ together with a measurable
real-valued function on $X$.

Any operator $A \in \Bc(H)$ can be decomposed into its ``real'' and
``imaginary'' parts, $A = {\rm Re}(A) + i{\rm Im}(A)$, where
${\rm Re}(A) = (A + A^*)/2$ and ${\rm Im}(A) = (A - A^*)/2i$ are both
self-adjoint, and this is analogous
to the decomposition of a complex-valued function into its real and
imaginary parts. Notice that in the finite-dimensional case, when $A$
is identified with a matrix, self-adjointness of $A$ does {\it not}
correspond to its entries being real. If we had used real scalars
from the start, we would still have to deal with non self-adjoint
operators and the general decomposition into real and imaginary
parts would not exist. This shows the intrinsically complex (as
opposed to real) character of Hilbert space operators and should
help explain why it is natural for us to use complex scalars.
(To emphasize this point further, note that there is also a ``polar
decomposition'' of arbitrary bounded operators which is analogous
to expressing a complex number in the form $z = re^{i\theta}$.)

A C*-algebra is a subset of $\Bc(H)$ that is closed under all
relevant structure.

\begin{defi}
A {\it concrete C*-algebra} is a linear subspace $\Ac$ of $\Bc(H)$
(for some Hilbert space $H$) that is closed in the norm topology and
satisfies
$$A, B \in \Ac \qquad \Rightarrow \qquad AB, A^* \in \Ac.$$
An {\it abstract C*-algebra} is a Banach space equipped with a product
and an involution that is linearly isometric to some concrete
C*-algebra, such that the isomorphism respects multiplicative and
involutive as well as linear structure.
\end{defi}

There is also an elegant abstract Banach algebra characterization
of C*-algebras. It can be used to show that the quotient of a
C*-algebra by any closed two-sided ideal is again a C*-algebra,
which is not obvious from the concrete characterization.
(Any such ideal is automatically stable under involution.)

(We shall also make occasional reference to von Neumann algebras.
Concretely, a {\it von Neumann algebra} is a C*-algebra that is closed
in the weak operator topology, the weakest topology on $\Bc(H)$ which
makes the function $A \mapsto \langle Av, w\rangle$ continuous for
every $v,w \in H$. It is {\it separably acting} if $H$ is separable.
Abstractly, a von Neumann algebra is a C*-algebra that is linearly
isometric to the dual of some Banach space, and the abstract
version of ``separably acting'' is that its predual is separable.)

C*-algebras play a basic role in mathematical physics, especially
quantum field theory \cite{Ara, BR, Wal}, and they have also had
major applications to other areas of mathematics, including index
theory for foliations \cite{Cnn}, the theory of group representations
\cite{Dix}, knot theory \cite{Jon}, and the Novikov conjecture in
algebraic topology \cite{Kas}.

\subsection{Basic examples}

\begin{exam}\label{B(H)}
$\Bc(H)$ is itself a C*-algebra. If $H = \Cb^n$ then we can
identify $\Bc(H)$ with $M_n(\Cb)$, the algebra of $n \times n$
complex matrices.

Direct sums of C*-algebras are again C*-algebras. Thus
$M_{n_1}(\Cb) \oplus \cdots \oplus M_{n_k}(\Cb)$ --- the set of
block diagonal matrices with blocks of size $n_1, \ldots, n_k$ ---
is a C*-algebra for any $k \in \omega$ and $n_1, \ldots, n_k \in
\omega$. Theorem: every finite-dimensional C*-algebra is isomorphic
to one of this form (\cite{W1}, Theorem 11.1.2).
\end{exam}

\begin{exam}\label{l^infty}
$l^\infty$, the space of all bounded sequences of
complex numbers, is a C*-algebra. Sums and products of
sequences are defined pointwise, the involution is pointwise
complex conjugation, and the norm of a sequence is its supremum
norm. There is an isometric embedding of $l^\infty$ into
$\Bc(l^2)$ that realizes elements of $l^\infty$ as
multiplication operators on $l^2$. Namely, $a = (a_i) \in l^\infty$ acts
on $v = (v_i) \in l^2$ by $av = (a_iv_i)$. This embedding preserves both
the norm and all algebraic structure, so $l^\infty$ is a C*-algebra.
I will call it the {\it diagonal embedding} because if we represent
operators on $l^2$ by infinite matrices, then the multiplication
operators just mentioned appear as diagonal matrices.
\end{exam}

\begin{exam}\label{C(X)}
Let $X$ be a compact metric space. Then $C(X)$, the space of all
continuous complex-valued functions on $X$, is a C*-algebra. Just
as with $l^\infty$, its algebraic operations are defined pointwise
and its norm is the supremum norm. We can embed $C(X)$ into $l^\infty$
by fixing a dense sequence $(x_n)$ in $X$ and taking a function
$f \in C(X)$ to the sequence $(f(x_n))$. Since we have already
seen that $l^\infty$ is a C*-algebra, this shows that $C(X)$ is
also a C*-algebra.

More generally, if $X$ is any compact Hausdorff space then $C(X)$
is a C*-algebra, but embedding it in $\Bc(H)$ might require that
$H$ be nonseparable. Theorem: every unital abelian C*-algebra is
isomorphic to one of this form (\cite{W1}, Theorem 5.3.5).
(A C*-algebra is {\it unital} if it contains a multiplicative unit.)
\end{exam}

Notice that $l^\infty$ is a unital abelian C*-algebra, so according
to the theorem just mentioned it must be isomorphic to $C(X)$ for
some compact Hausdorff space $X$. Indeed, $l^\infty \cong C(\beta\omega)$.

I mentioned earlier that the quotient of a C*-algebra by any closed
two-sided ideal is always a C*-algebra. For example, $c_0$, the
space of all complex sequences which converge to zero, is a closed
two-sided ideal of $l^\infty$, so the quotient space $l^\infty/c_0$
is also a (unital abelian) C*-algebra. Indeed, $l^\infty/c_0 \cong
C(\beta\omega - \omega)$.

\begin{exam}\label{C(l^2)}
An operator $A \in \Bc(H)$ is said to be of {\it finite rank} if
its range is finite-dimensional. The set of all finite rank operators
on $H$ is stable under sums, products, and adjoints, but if $H$
is infinite-dimensional it fails to be a C*-algebra because it
is not closed in norm. However, its closure, denoted $\Kc(H)$, is a
C*-algebra. The operators in $\Kc(H)$ are precisely those bounded
operators that can be approximated in norm by finite rank operators;
they are called {\it compact} operators because they are also
characterized by the fact that the image of the unit ball of $H$ under
such an operator has compact closure.

$\Kc(H)$ is a C*-algebra and it is also a closed two-sided ideal
of $\Bc(H)$. When $H = l^2$ the quotient C*-algebra $\Cc(l^2) =
\Bc(l^2)/\Kc(l^2)$ is called the {\it Calkin algebra}.
\end{exam}

Recall the diagonal embedding of $l^\infty$ into the bounded operators
on $l^2$ (Example \ref{l^infty}). The diagonal operators which are compact
are precisely those whose diagonal entries converge to zero, and hence
they correspond to elements of $c_0 \subseteq l^\infty$. In other
words, identifying $l^\infty$ with its image in $\Bc(l^2)$, we have
$l^\infty \cap \Kc(l^2) = c_0$.

I suggested earlier that self-adjoint operators
are analogous to real-valued functions. This is one
link in an extended chain of analogies which compare closed
subsets of a set with closed subspaces of a Hilbert space,
spaces of continuous functions on compact Hausdorff spaces
with C*-algebras, $L^\infty$ spaces with von Neumann algebras,
etc. This is all motivated at a fundamental level by the
transition from classical to quantum mechanics, so that one
sometimes hears C*-algebras referred to as ``quantum topological
spaces'', and so on. A readable account of this motivation is
given in Chapter 1 of \cite{W1}.

In particular, $\Bc(l^2)$ can be regarded as a
noncommutative analog of $l^\infty \cong C(\beta\omega)$, $\Kc(l^2)$
as an analog of $c_0 = C_0(\omega)$ (the continuous functions on
$\omega$ that vanish at infinity), and $\Cc(l^2)$ as an analog of
$l^\infty/c_0 \cong C(\beta\omega-\omega)$. Because of the
correspondence between topological spaces and abelian C*-algebras
one might simply say that $\Bc(l^2)$, $\Kc(l^2)$, and $\Cc(l^2)$
are analogs of $\beta\omega$, $\omega$, and $\beta\omega - \omega$.
I will return to this point in Section 3.

\begin{exam}\label{Fermion}
For each $n \in \omega$ define an embedding $\phi_n: M_{2^n}(\Cb)
\to M_{2^{n+1}}(\Cb)$ by
$$\phi_n(A) = \left[
\begin{matrix}
A&0\cr
0&A
\end{matrix}
\right].$$
This is a linear isometry from $M_{2^n}(\Cb)$ into $M_{2^{n+1}}(\Cb)$
which respects both multiplicative and involutive structure, so with
a slight abuse of notation we can define $\Fc_0 =
\bigcup_1^\infty M_{2^n}(\Cb)$. Taking the Cauchy completion then
yields a C*-algebra $\Fc$, the {\it Fermion} or {\it CAR} algebra.
(``CAR'' stands for ``canonical anticommutation relations''.)
Concretely, $\Fc$ can be realized as acting on a Hilbert space
$H$ that has an orthonormal basis indexed by the finite subsets
of $\omega$. To see this, observe that the standard basis of
$\Cb^{2^n}$ can be indexed by the subsets of $n$.
Thus $H$ naturally decomposes as a direct sum of copies of $\Cb^{2^n}$,
one for each finite subset of $\omega - n$, and we
can let $M_{2^n}(\Cb)$ act on each copy. These actions match up with
the embeddings $\phi_n$, so $\Fc$ can be defined to be the closure of
the union over $n \in \omega$ of the images in $\Bc(H)$ of $M_{2^n}(\Cb)$.
\end{exam}

The Fermion algebra might be considered an analog of (the
continuous functions on) the Cantor
set $K$. The intuition here comes from associating a diagonal
$2^n \times 2^n$ matrix $A = {\rm diag}(a_1, \ldots, a_{2^n})$
with the function $f \in C(K)$ which takes the value $a_i$ on
the $i$th subinterval of length $3^{-n}$. In this way we obtain
a natural isomorphism of the ``diagonal'' of $\Fc$, i.e., the
closure of the union over $n \in \omega$ of the diagonal matrices
in $M_{2^n}(\Cb) \subset \Fc$, with the continuous functions on $K$.

\subsection{States and representations}

Although in practice most C*-algebras have some simple favored
realization in some special $\Bc(H)$ (the Calkin algebra is a
notable exception), the theory of representations of C*-algebras
is nonetheless of great intrinsic interest and of fundamental
importance, for example, in applications of C*-algebras to physics
\cite{BR}, group representations \cite{Dix}, and wavelets \cite{Jor}.

\begin{defi}
Let $\Ac$ be a C*-algebra and let $H$ be a Hilbert space. A
{\it representation} of $\Ac$ on $H$ is a linear map
$\phi: \Ac \to \Bc(H)$ which respects the product and involution
on $\Ac$. It is {\it unital} if $\Ac$ is unital and
$\phi$ takes the unit of $\Ac$ to the identity operator on $H$.
It is {\it reducible} if $H$ can be decomposed as a direct sum of
two orthogonal proper subspaces each of which is invariant
for the action of $\Ac$ (i.e., $\phi(A)$ takes each of the two subspaces
into itself, for every $A \in \Ac$); otherwise it is {\it irreducible}.
\end{defi}

Every representation is automatically continuous (in fact, contractive).

Irreducible representations are the ``building blocks'' of representation
theory, and morally we think of all representations as being in some sense
decomposable
into irreducible representations, although this sentiment can be made
precise only for separable C*-algebras and even then the exact statement is
somewhat complicated. Also, if one approaches C*-algebras from the perspective
of their being a noncommutative generalization of the spaces $C(X)$,
then the irreducible representations are for some purposes
a good generalization of the notion ``points of the underlying space $X$''.
This is because the irreducible representations of $C(X)$, up to the
natural notion of equivalence, canonically correspond to the points of $X$.
Namely, for each $x \in X$ we have an evaluation functional
$\phi: C(X) \to \Cb \cong M_1(\Cb)$ given by $\phi(f) = f(x)$, and this is
an irreducible representation of $C(X)$ on a one-dimensional Hilbert space,
and every irreducible representation of $C(X)$ is equivalent to one of
these. (Two representations $\phi_i:
\Ac \to \Bc(H_i)$, $i = 1,2$, are {\it equivalent} if there is a Hilbert
space isomorphism $U: H_1 \cong H_2$ such that $\phi_1(A) = U^{-1}\phi_2(A)U$
for all $A \in \Ac$.)

A very useful tool for dealing with representations is provided
by the concept of a state on a C*-algebra. Before explaining this
notion we first need to define positivity. Let $\Ac \subseteq \Bc(H)$
be a C*-algebra and let $A \in \Ac$. Then we say that $A$ is
{\it positive}, and we write $A \geq 0$, if any of the following
equivalent conditions is satisfied: (1) $\langle Av,v\rangle \geq 0$
for all $v \in H$; (2) $A = B^*B$ for some $B \in \Ac$; (3) there is
a measure space $(X,\mu)$, a bounded measurable function $f: X \to
[0,\infty)$, and a Hilbert space isomorphism from $H$ to $L^2(X,\mu)$
taking $A$ to $M_f$. Further, for $A, B \in \Ac$ self-adjoint we
set $A \leq B$ if $B - A$ is positive. This is a partial order on
the self-adjoint elements of $\Ac$.

\begin{defi}
Let $\Ac$ be a unital C*-algebra. A linear functional $\rho: \Ac \to \Cb$
is a {\it state} if (1) $\rho(I_\Ac) = 1$, where $I_\Ac$ is the unit of
$\Ac$; (2) $A \geq 0$ implies $\rho(A) \geq 0$; and (3) $\|\rho\| = 1$,
where $\|\rho\|$ is defined to be $\sup\{|\rho(A)|/\|A\|: A \neq 0\}$.
It is {\it pure} if it is an extreme point of the set of states, i.e.,
it cannot be expressed as a convex combination of two other states.
\end{defi}

Actually any two of the conditions in the definition of a state imply
the third. For simplicity we omit discussion of states on non-unital
C*-algebras, although generalization to that case is reasonably
straightforward. (For example, every non-unital C*-algebra can be
trivially embedded in a unital C*-algebra, so that one can effectively
reduce to the unital case.)

The states on $C(X)$ correspond to Borel probability measures on $X$.
For any such measure $\mu$, the map $f \mapsto \int f\, d\mu$ is a
state on $C(X)$, and all states arise in this way. The pure states
are the evaluation functionals, which correspond to probability
measures concentrated at a single point.

The prototypical example of a state on a general unital C*-algebra
arises in the following way. Let
$\phi: \Ac \to \Bc(H)$ be a unital representation and let $v \in H$ be
a unit vector. Then the map $\rho: \Ac \to \Cb$ defined by $\rho: A \mapsto
\langle \phi(A)v,v\rangle$ is a state on $\Ac$, called a {\it vector state}
(for that representation). Thus, representations give rise to states.

Conversely, states give rise to representations via the {\it GNS
(Gelfand-Naimark-Segal) construction}, executed as follows. Let $\Ac$
be a unital
C*-algebra and let $\rho$ be a state on $\Ac$. For any $A,B \in \Ac$
define $\langle A,B\rangle = \rho(B^*A)$. In general this is only a ``pseudo
inner product'' because nonzero vectors might have norm zero. But by
factoring out null vectors and completing we obtain a Hilbert space
$H_\rho$, and $\Ac$ acts on this Hilbert space by left multiplication.
The resulting representation is denoted $\phi_\rho: \Ac \to \Bc(H_\rho)$.

Let $v \in H_\rho$ be the (unit) vector that comes from the unit of
$\Ac$. Then $\langle \phi_\rho(A)v, v\rangle = \rho(A)$ for all $A \in \Ac$.
Thus, we recover the original state $\rho$ as a vector state.
So we can pass back and forth between states and representations.
Theorem: The state $\rho$ is pure if and only if the representation
$\phi_\rho$ is irreducible, and every irreducible representation is
equivalent to $\phi_\rho$ for some pure state $\rho$ (\cite{Ped},
Theorem 3.13.12).

An operator $U \in \Bc(H)$ is {\it unitary} if $U^*U = UU^* = I$;
equivalently, $U: H \to H$ is a Hilbert space isomorphism, i.e., a
bijection that preserves the inner
product. We say that two pure states $\rho_1$ and $\rho_2$ on a unital
C*-algebra $\Ac$ are {\it equivalent} if there is a unitary $U \in \Ac$
such that $\rho_1(A) = \rho_2(U^*AU)$ for all $A \in \Ac$. Theorem:
$\rho_1$ and $\rho_2$ are equivalent pure states if and only if
$\phi_{\rho_1}$ and $\phi_{\rho_2}$ are equivalent irreducible
representations (\cite{Ped}, Proposition 3.13.4).

Let $\Ac$ be a unital C*-algebra and let $\Bc \subseteq \Ac$ be a
C*-algebra contained in $\Ac$ that contains the unit of $\Ac$. Then
any state on $\Ac$ restricts to a state on $\Bc$ (this is obvious),
and by the Hahn-Banach theorem any state on $\Bc$ extends to at least
one state on $\Ac$. Moreover, using convexity theory it is not too hard
to show that any pure state on $\Bc$ can be extended to at least one
pure state on $\Ac$. However, a pure state on $\Ac$ might not restrict
to a pure state on $\Bc$. For example, let $\Ac = \Bc(L^2[0,1])$ and let
$\Bc = C[0,1]$, embedded in $\Ac$ as
multiplication operators. Then for any unit vector $\psi \in L^2[0,1]$
the corresponding vector state on $\Ac$ is pure since every vector
state on any $\Bc(H)$ is pure, but its restriction to $\Bc$ is the
state $f \mapsto \int f|\psi|^2\, dx$, which is not pure.

\section{Consistency results for C*-algebras}

I will now outline a number of consistency results appearing in or near
the C*-algebra literature. In some cases they are explicitly formulated
as consistency statements, but often they are simply presented as
theorems proven with the aid of CH.

Exactly what constitutes ``the C*-algebra literature'' is debatable.
The first three topics discussed below are slightly questionable in
this respect, though all three should be intelligible and of interest
to C*-algebraists. The last three topics in this section deal with von
Neumann algebras.

\subsection{Automatic continuity}

Probably the best known independence result in functional analysis
states that if $X$ is an infinite compact Hausdorff space then it is
independent of ZFC whether there exists a discontinuous homomorphism
(= linear multiplicative map) from $C(X)$ into a Banach algebra \cite{DW}.
Although $C(X)$ is a C*-algebra, this is more of a Banach algebra result
than a C*-algebra result because the involution on $C(X)$ is not in play
and the range algebra is well out of the realm of C*-algebras.

There is a substantial literature on automatic continuity, some of
which specifically considers noncommutative C*-algebras (being mapped
into Banach algebras). Assuming CH, Laursen \cite{Lau} showed
that a separable C*-algebra has a discontinuous homomorphism into a
commutative Banach algebra if and only if the (non-closed) two-sided
ideal generated by all commutators has infinite codimension.
Runde \cite{Run} used CH to construct a discontinuous homomorphism
from a noncommutative C*-algebra which has a pathological property,
namely, that it fails to decompose in a way that is possible for any
discontinuous homomorphism from a commutative Banach algebra. 

\subsection{The second dual of $C[0,1]$}

The Banach space dual of $C[0,1]$ is $M[0,1]$, the space of
complex Borel measures on $[0,1]$. The dual of $M[0,1]$ is the
``enveloping von Neumann algebra'' of $C[0,1]$ and is fairly
intractable. However, assuming CH, R.\ D.\ Mauldin \cite{Mau} proved
the following representation theorem for bounded linear functionals
on $M[0,1]$: for every such functional $T$ there is a bounded
function $\psi$ from the set of Borel subsets of $[0,1]$ to $\Cb$
such that
$$T(\mu) = \int \psi\, d\mu$$
for all $\mu \in M[0,1]$. Here the integral notation signifies
the limit over the directed set of Borel partitions $\{B_1, \ldots,
B_n\}$ of $[0,1]$ of the quantity $\sum \psi(B_i)\mu(B_i)$, and
part of the assertion is that $\psi$ can be chosen so that this
limit always exists. The proof begins by choosing a maximal family
of mutually singular Borel measures on $[0,1]$; using CH, this
family is then indexed by countable ordinals and $\psi$ is then
defined by transfinite recursion.

\subsection{Sum decomposition of operator ideals}

The compact operators on $l^2$ constitute a two-sided ideal of
$\Bc(l^2)$, and this is
the only nontrivial closed two-sided ideal. However, there are many
nontrivial non-closed ideals, each contained in $\Kc(l^2)$ and containing
all finite rank operators. Of special interest are the {\it Schatten
$p$-ideals} $\Sc_p(l^2)$ defined for $1 \leq p < \infty$ as the set of
compact operators the eigenvalues of whose real and imaginary parts
are both $p$-summable; these are regarded as Hilbert space analogs
of the Banach spaces $l^p$.

The (not necessarily closed) ideals of $\Bc(l^2)$ form a complete
lattice under inclusion. Blass and Weiss
\cite{BW} proved that CH implies that every proper ideal that properly
contains the finite rank operators is the sum (join) of two smaller
ideals. They also noted that MA is actually sufficient for their proof.
Subsequently Blass and Shelah \cite{Bla, BS} established
that the failure of $\Kc(l^2)$ to be the join of two smaller ideals
is equivalent to the assertion of near coherence of ultrafilters
(NCF), i.e., any two free ultrafilters over
$\omega$ are related (in both directions) by a finite-to-one map
from $\omega$ to itself, and that the latter condition is consistent
with ZFC.

\subsection{Ultrapowers of C*-algebras}

Ultraproduct constructions are useful for studying approximation
properties in C*-algebras and have been heavily used, especially
in von Neumann algebra theory. This motivated Ge and Hadwin to
initiate a general study of C*-algebra ultraproducts \cite{GH},
which are a straightforwardly obtained from set-theoretic
ultraproducts by removing elements of infinite norm and factoring
out elements of norm zero.
They showed that CH is equivalent to the assertion that
all nontrivial ultrapowers of any separable C*-algebra are isomorphic.
In slightly more detail, CH implies that if $\Ac$ is a separable
C*-algebra, then any two ultrapowers of $\Ac$ with respect to free
ultrafilters over $\omega$ are isomorphic by a map which fixes
the elements of $\Ac$; conversely, $\neg$CH implies that any separable
C*-algebra whose set of projections satisfies a simple condition has
non-isomorphic ultrapowers with respect to free ultrafilters over $\omega$.
The second statement follows easily from a theorem of Dow and Shelah
\cite{Dow, S2} that $\neg$CH implies that any partially ordered set
containing an infinite chain has nontrivial ultrapowers which are not
order isomorphic.

A related result of Ge and Hadwin states that CH implies that all
nontrivial tracial ultrapowers (a variation on the C*-algebra
ultrapower construction) of a separably acting finite von Neumann
algebra are isomorphic. This answered a thirty year old
question of McDuff.

Ge and Hadwin also proved that the existence of a measurable cardinal
is equivalent to the existence of a nontrivial separable ultrapower
of a separable C*-algebra.

\subsection{Naimark's problem}

Naimark \cite{N1} showed that every irreducible representation of
$\Kc(H)$ is equivalent to the identity representation on $H$,
and he asked \cite{N2} whether there are any other C*-algebras with
only one irreducible representation up to equivalence. Rosenberg
\cite{Ros} soon showed that there are no other separable C*-algebras
with this property. In subsequent work of Glimm \cite{Gli} and others
this result was subsumed in a general theory of representations of
separable C*-algebras. General C*-algebras were classified into two
categories, type I and not type I; the irreducible representations
of type I C*-algebras were in some sense completely described,
whereas separable non-type I C*-algebras were shown to have a
uniformly bad representation theory. (The simplest definition of
``type I'' is: every irreducible representation on a Hilbert space
$H$ contains $\Kc(H)$. However, there are many equivalent conditions
framed in very different terms; see Theorem 6.8.7 of \cite{Ped}.) In
particular, Glimm showed that any non-type I C*-algebra contains a
C*-algebra that has the Fermion algebra $\Fc$ (Example \ref{Fermion})
as a quotient. Since $\Fc$ has $2^{\aleph_0}$ inequivalent
irreducible representations it easily follows that the same is
true of any separable non-type I C*-algebra.

Naimark's question for nonseparable C*-algebras remained open for
fifty years. This seemed odd since it was known that there could be no
type I counterexamples, whereas one might assume that the space of
representations of any nonseparable non-type I C*-algebra should be,
if anything, even less tractable as a topological space than the
representation spaces of separable non-type I C*-algebras. (The set
of all representations up to equivalence carries a natural topology.)
However, this expectation is naive because even though
nonseparable non-type I C*-algebras ought to have ``many'' pure states,
they might also contain ``many'' unitaries making these pure states
equivalent. There is therefore an evident possibility of constructing
a nonseparable counterexample to Naimark's problem by, loosely speaking,
repeatedly adding unitaries until all pure states become equivalent.

This was done in \cite{AW} using diamond to get around the problem
that as unitaries are added the set of pure states grows.  The
key technical step is a construction which enables one to add a
unitary to any simple, separable, unital C*-algebra in such a way
that two chosen pure states extend uniquely to the larger
algebra and are made equivalent. A recursive construction of
length $\omega_1$ can then be carried out where at each step
some new pure state, selected using diamond, is made equivalent
to a pure state that is fixed throughout the construction. Ultimately
all pure states are equivalent since any pure state on the resulting
C*-algebra restricts to a pure state on a closed unbounded family of
C*-algebras appearing in the construction; diamond then ensures that
every pure state will have been made equivalent to the fixed
pure state at some stage.

By Glimm's result mentioned earlier, any counterexample to Naimark's
problem cannot be generated by fewer than $2^{\aleph_0}$ elements.
As the counterexample constructed using diamond is generated by
$\aleph_1$ elements, it follows that the existence of an
$\aleph_1$-generated counterexample is independent of ZFC. However,
it is presumably the case that the nonexistence of any counterexample
is consistent with ZFC. This question remains open. Even better would
be the consistency of every nonseparable non-type I C*-algebra having
uncountably many inequivalent irreducible representations.

\subsection{Pure states on von Neumann algebras}

Akemann, Anderson, and Pedersen \cite{AAP1, AAP2} used CH to
investigate the behavior of pure states on von Neumann algebras.
In \cite{AAP1} they showed that any sequence of pure states on a
separably acting von Neumann algebra that is supported by a sequence
of mutually orthogonal, positive, norm-one elements has at least one
limit point which is a pure state. This was done by using CH to
construct an ultrafilter with special properties in order
to identify a limiting pure state. The authors remark that their
proof ``could proceed with an axiom which is strictly weaker than
the Continuum Hypothesis,'' probably a reference to MA. What is
actually needed for the proof is the nonexistence of an inextendible
$\kappa$-tower in $P(\omega)/{\rm fin}$ for any $\kappa < 2^{\aleph_0}$.

The result in \cite{AAP2} is related to the concept of a
{\it perfect} C*-algebra \cite{AS}. It is shown in \cite{AAP2}
that a separable C*-algebra is perfect if and only if it has
only trivial diffuse sequences, where a sequence $(A_n)$ in a
unital C*-algebra $\Ac$ is {\it diffuse} if every net of pure
states $(\rho_\alpha)$ on $\Ac$ that converges weak* to a pure
state satisfies $\lim_{\alpha, n} \rho_\alpha(A_n^*A_n + A_nA_n^*) = 0$,
and a diffuse sequence is {\it trivial} if $A_n \to 0$. It is easier
for nonseparable C*-algebras to be perfect \cite{Arc}. The authors
conjecture that every von Neumann algebra is perfect, and as evidence
for this conjecture prove that assuming CH, every separably acting
von Neumann algebra has only trivial diffuse sequences. This result
is proven by first showing that the existence of a nontrivial diffuse
sequence implies the existence of a mutually orthogonal, positive,
norm one diffuse sequence, and then invoking the theorem in \cite{AAP1}
described above to reach a contradiction.

\subsection{Quotients of finite von Neumann algebras}

A von Neumann algebra $\Mx \subseteq \Bc(H)$ is {\it finite} if it does
not contain an isometry from $H$ onto a proper subspace of $H$. It is
{\it properly infinite} if $\Mx \cong \Mx_1 \oplus \Mx_2$ with $\Mx_1$
finite implies $\Mx_1 = \{0\}$. Every von Neumann algebra
is the direct sum of a finite von Neumann algebra and a properly
infinite von Neumann algebra.

J.\ Vesterstr{\o}m \cite{Ves} investigated the question of when the
quotient of a finite von Neumann algebra by a norm closed two-sided
ideal is a von Neumann algebra. The problem had previously been
answered by Takemoto for separably acting properly infinite
von Neumann algebras; in that case the quotient is a von Neumann
algebra if and only if the ideal is weak operator closed. The finite
case is more complicated, but assuming CH Vesterstr{\o}m gave reasonable
necessary and sufficient conditions for the quotient of a separably
acting finite von Neumann algebra by a norm closed two-sided ideal
to be a von Neumann algebra. He used CH in proving the necessity
of the condition that the quotient of the center by its intersection
with the ideal in question must be of the form $L^\infty(Y,\nu)$ with $(Y,\nu)$
a $\sigma$-finite measure space. This comes down to showing that if
$L^\infty(Y,\nu)$ is a quotient of $L^\infty(X,\mu)$, $L^\infty(X,\mu)
\cong L^1(X,\mu)^*$ and $L^\infty(Y,\nu) \cong L^1(Y,\nu)^*$, and $L^1(X,\mu)$
is separable, then $(Y,\nu)$ is $\sigma$-finite, which follows from
CH by a simple cardinality argument. Actually, all that is needed for
this argument is $2^{\aleph_0} < 2^{\aleph_1}$.

\subsection{Weak expectations}

Blackadar \cite{Blk} used CH to construct a noninjective separably
acting factor von Neumann algebra which has a weak expectation.
The technique involves the construction of a nuclear C*-algebra
by a transfinite recursion of length $\omega_1$, with CH being invoked
to ensure that there are only $\aleph_1$ simple separable C*-algebras
to be handled by the construction. See \cite{Blk} for explanation of
the terms used here. As Blackadar notes, it is likely that the
assumption of CH can be dropped.

\section{The Calkin algebra}

Recall from Section 1 that the Calkin algebra is the quotient of
the bounded operators on $l^2$ by the compact operators,
$\Cc(l^2) = \Bc(l^2)/\Kc(l^2)$. It will become apparent below that
$\Cc(l^2)$ can fruitfully be thought of as a noncommutative analog
of (the continuous functions on) $\beta\omega - \omega$. For example,
an old open problem asks whether every automorphism of $\Cc(l^2)$
is inner, i.e., implemented by a unitary $u \in \Cc(l^2)$. In more
detail: let $\phi: \Cc(l^2) \to \Cc(l^2)$ be a linear bijection
which respects the product and involution; must there exist
$u \in \Cc(l^2)$ such that $\phi(a) = u^*au$ for all $a \in \Cc(l^2)$?
This problem is a natural analog of the question whether there exist
nontrivial self-homeomorphisms of $\beta\omega -\omega$, which suggests
that it may be independent of ZFC and that the techniques developed for
proving independence of the classical question \cite{S1, SS, vM, Vel}
would be useful in showing this.

\subsection{Projections in $\Cc(l^2)$}

For any closed subspace of a Hilbert space $H$ there is a bounded operator
on $H$ which orthogonally projects vectors onto the subspace. Such
an  operator is called an {\it orthogonal projection}, or often simply
a {\it projection}. Projections are algebraically characterized by
the fact that they are self-adjoint idempotents, i.e., they
satisfy $P = P^2 = P^*$.

In a unital abelian C*-algebra $C(X) \subset \Bc(H)$, this algebraic
characterization shows that the projections are identified with the
continuous functions on $X$ whose range is contained in $\{0,1\}$. These
are precisely the characteristic functions of clopen subsets of
$X$. Thus the set of projections in $C(X)$ is a Boolean algebra
that is naturally isomorphic to the Boolean algebra of clopen
subsets of $X$. If $X$ is a Stone space (i.e., is totally disconnected
compact Hausdorff) then the linear span of the set of projections
in $C(X)$ is dense in $C(X)$, and working with the Boolean algebra
of clopen subsets of $X$ is more or less equivalent to working with
the C*-algebra $C(X)$. For general compact Hausdorff spaces, Stone
duality with Boolean algebras is therefore, broadly speaking,
generalized by Gelfand duality with unital abelian C*-algebras;
see, e.g., Section 5.1 of \cite{W1}.

Now $\beta\omega$ and $\beta\omega - \omega$ are Stone spaces, and
their dual Boolean algebras are $P(\omega)$ and $P(\omega)/{\rm fin}$.
Thus, we may regard the lattices of projections in $\Bc(l^2)$ and
$\Cc(l^2)$, which we denote $\Pc(\Bc(l^2))$ and $\Pc(\Cc(l^2))$, as
Hilbert space analogs of $P(\omega)$ and $P(\omega)/{\rm fin}$.

There is a fairly substantial literature on chains of projections in
$\Bc(l^2)$ (they are called ``nests''; see \cite{Dav}). The study of
chains of projections in $\Cc(l^2)$ was initiated by Hadwin in \cite{Had},
where he proved that CH implies that all maximal chains are order-isomorphic.
The proof goes by showing that any maximal chain is $\aleph_0$-saturated
and then invoking a classical result of Hausdorff. Hadwin conjectured
that isomorphism of all maximal chains of projections in $\Cc(l^2)$ is
equivalent to CH, which seems unlikely. My student Eric Wofsey \cite{Wof}
has recently shown that the standard notion of forcing used to force
inextendible towers in $P(\omega)/{\rm fin}$ also forces inextendible
towers in $\Pc(\Cc(l^2))$, and he used this to prove the consistency
of the existence of non-isomorphic maximal chains in $\Pc(\Cc(l^2))$.
Wofsey is also studying maximal
families of orthogonal projections in $\Cc(l^2)$, i.e., families of
projections the product of any two of which is zero. He calls such
a family a {\it maximal almost orthogonal family} or {\it maof},
and he has shown that it is consistent with ZFC that there is an
maof of cardinality strictly less than $2^{\aleph_0}$.

In studying the lattice of projections in $\Cc(l^2)$ it is possible
to largely bypass the machinery of C*-algebras. That is because this
lattice can be identified with the lattice of projections in
$\Bc(l^2)$, i.e., the closed subspaces of $l^2$, modulo a natural
equivalence relation. For the reader's benefit we prove this folklore
result. We use the following easily verified fact: suppose $(X,\mu)$
is a measure space and $f \in L^\infty(X,\mu)$, and let
$M_f \in \Bc(L^2(X,\mu))$ be the operator of multiplication by $f$.
Then $M_f$ is compact if and only if for every $\epsilon > 0$ the
space $L^2(X_\epsilon, \mu|_{X_\epsilon})$ is finite-dimensional,
where $X_\epsilon = \{x \in X: |f(x)| \geq \epsilon\}$.

Let $\pi: \Bc(l^2) \to \Cc(l^2)$ be the natural quotient map.

\begin{prop}
Let $p \in \Cc(l^2)$ be a projection. Then there exists a projection
$P \in \Bc(l^2)$ such that $\pi(P) = p$.
\end{prop}

\begin{proof}
We have $p = p^* = p^2$. Choose $A \in \Bc(l^2)$ such that $\pi(A) = p$. Then
$$\pi(A^*) = \pi(A)^* = p^* = p$$
so $\pi(B) = p$ where $B = (A + A^*)/2$. Now $B$ is self-adjoint so
we can replace $\Bc(l^2)$ by $\Bc(L^2(X,\mu))$ for some measure space
$(X,\mu)$ and $B$ by $M_f$ for some real-valued $f \in L^\infty(X,\mu)$.

We have $\pi(B) = p = p^2 = \pi(B^2)$, so $B - B^2$ is compact. That
is, $M_g$ is a compact operator where $g = f - f^2$. Define
$$h(x) = 
\begin{cases}
1&\mbox{if $f(x) \geq 1/2$}\\
0&\mbox{if $f(x) < 1/2$}.
\end{cases}$$
We claim that $M_f - M_h = M_{f - h}$ is compact, so that $P = M_h$ is the
desired projection such that $\pi(P) = p$. Thus, fix
$\epsilon > 0$ and let $X_\epsilon = \{x \in X: |f(x) - h(x)|
\geq \epsilon\}$; we must show that $L^2(X_\epsilon, \mu|_{X_\epsilon})$
is finite-dimensional. But assuming $\epsilon \leq 1/2$, we have
$$X_\epsilon \subseteq \{x \in X: |g(x)| \geq \epsilon - \epsilon^2\}$$
and therefore compactness of $M_g$ implies that
$L^2(X_\epsilon,\mu|_{X_\epsilon})$ is finite-dimensional.
\end{proof}

(For more results along these lines see Theorem 4.3 of \cite{AP} and
the references cited there.)

\begin{defi}
Let $E$ and $F$ be closed subspaces of
$l^2$. Say that $E$ is {\it essentially contained in $F$} if for every
$\epsilon > 0$ there exists a cofinite-dimensional subspace $E_0$ of
$E$ such that every unit vector $v$ in $E_0$ satisfies $d(v, F) \leq
\epsilon$.
\end{defi}

Here $d(v,F)$ is the distance between the vector $v$ and the subspace
$F$. Equivalently, we could ask that for every $\epsilon > 0$ there
is a cofinite-dimensional subspace $E_0$ of $E$ such that $\|P_Fv\|
\geq 1 - \epsilon$ for every unit vector $v \in E_0$, where $P_F$ is
the orthogonal projection onto $F$.

\begin{prop}
Let $P$ and $Q$ be projections in $\Bc(l^2)$. Then $\pi(P) \leq \pi(Q)$
if and only if ${\rm ran}(P)$ is essentially contained in ${\rm ran}(Q)$.
\end{prop}

\begin{proof}
($\Rightarrow$) Suppose $E = {\rm ran}(P)$ is not essentially contained in
$F = {\rm ran}(Q)$. Then
we can find $\epsilon > 0$ and an orthogonal sequence of unit vectors
$(v_n)$ in $E$ such that $d(v_n, F) \geq \epsilon$ for all $n$. If
$\pi(P) \leq \pi(Q)$ then $r = \pi(Q) - \pi(P)$ is a projection and
there exists a projection $R \in \Bc(l^2)$ such that $\pi(R) = r$. Then
$\pi(P) + \pi(R) = \pi(Q)$ so $P + R - Q$ is compact. But
\begin{eqnarray*}
\langle(P + R - Q)(v_n), v_n\rangle
&=& \langle v_n, v_n\rangle + \langle Rv_n, v_n\rangle
- \langle Qv_n, v_n\rangle\\
&\geq& 1 + 0  - (1 - \epsilon^2) = \epsilon^2,
\end{eqnarray*}
so we cannot have $(P + R - Q)(v_n) \to 0$. This contradicts
compactness of $P + R - Q$, so we must have $\pi(P) \not\leq \pi(Q)$.

($\Leftarrow$) Suppose $E$ is essentially contained in $F$.
Let $E'$ be the closure of $Q[E] \subseteq F$
and let $P'$ be the orthogonal projection onto $E'$. We claim that
$P - P'$ is compact, i.e., $\pi(P) = \pi(P')$; since $P' \leq Q$
this will imply $\pi(P) \leq \pi(Q)$.

(Note that $\pi$ preserves order because $A \geq 0$ $\Rightarrow$ $A = B^*B$
for some $B$ $\Rightarrow$ $\pi(A) = \pi(B^*B) = \pi(B)^*\pi(B) \geq 0$.)

Let $\epsilon > 0$. It will suffice to find a finite rank operator $A$ such
that $\|P - P' - A\| \leq \delta = \epsilon + 2\sqrt{\epsilon}$. To do
this let $E_0$ be a
cofinite-dimensional subspace of $E$ such that $d(v,F) \leq \epsilon$
for every unit vector $v$ in $E_0$. Let $P_0$ be the orthogonal
projection onto $E_0$ and let $P_0'$ be the orthogonal projection
onto $E_0' =$ the closure of
$Q[E_0] \subseteq F$. Then $A = (P - P_0) - (P' - P_0')$ is a
finite rank operator, so we must show that $\|P_0 - P_0'\| \leq \delta$.
Let $w \in l^2$ be any unit vector. Then $w_0 = P_0(w)$ is the closest
vector to $w$ in $E_0$ and similarly for $w_0' = P_0'(w)$. But there
exists $w_0'' \in E_0'$ such that $\|w_0 - w_0''\| \leq \epsilon$.
So
$$\|w - w_0'\| \leq \|w - w_0''\| \leq \|w - w_0\| + \epsilon,$$
and similarly $\|w - w_0\| \leq \|w - w_0'\| + \epsilon$. Since
$w - w_0'$ is perpendicular to $w_0' - w_0''$ we then have
\begin{eqnarray*}
\|w_0' - w_0''\|^2 &=& \|w - w_0''\|^2 - \|w - w_0'\|^2\\
&\leq& (\|w - w_0\| + \epsilon)^2 - (\|w - w_0\| - \epsilon)^2\\
&\leq& 4\epsilon.
\end{eqnarray*}
So
$$\|(P_0 - P_0')(w)\| = \|w_0 - w_0'\|
\leq \|w_0 - w_0''\| + \|w_0' - w_0''\|
\leq \epsilon + 2\sqrt{\epsilon} = \delta,$$
and as $w$ was arbitrary this shows that $\|P_0 - P_0'\| \leq \delta$.
\end{proof}

\begin{coro}
The lattice of projections in $\Cc(l^2)$ is isomorphic to the
lattice of closed subspaces of $l^2$ modulo equivalence, where
$E$ and $F$ are equivalent if each is essentially contained
in the other.
\end{coro}

It should be emphasized that the equivalence relation appearing
here is strictly weaker than the more obvious relation according
to which two closed subspaces of $l^2$ are equivalent if each is
contained in a finite-dimensional extension of the other. (For
example, letting $(e_n)$ be the standard basis of $l^2$, the closed
subspace generated by the vectors $e_{2n}$ is equivalent in the weaker
sense, but not the stronger sense, to the closed subspace generated by
the vectors $e_{2n} + \frac{1}{n}e_{2n+1}$.) The
latter also gives rise to an analog of $P(\omega)/{\rm fin}$, but
the C*-algebraic motivation for the former indicates that it is
of greater intrinsic interest.

\subsection{Masas and pure states on $\Cc(l^2)$}

A {\it masa} (``maximal abelian self-adjoint subalgebra'') of a C*-algebra
is a maximal abelian C*-subalgebra. In $\Bc(l^2)$ we have the following
characterization: a C*-algebra $\Ac \subseteq \Bc(l^2)$ is a masa if and
only if there is a $\sigma$-finite measure space $(X,\mu)$ and a
Hilbert space isomorphism $U: l^2 \cong L^2(X,\mu)$ taking $\Ac$ to
$L^\infty(X,\mu)$ (acting as multiplication operators on $L^2(X,\mu)$).
That is, $\Ac = \{U^{-1}M_fU: f \in L^\infty(X,\mu)\}$.

Johnson and Parrott \cite{JP} showed that if $\Ac$ is a masa of
$\Bc(l^2)$ then $\pi(\Ac)$ is a masa of $\Cc(l^2)$, where as before
$\pi: \Bc(l^2) \to \Cc(l^2)$ is the natural quotient map. In particular,
the canonical masa of $\Bc(l^2)$, namely $l^\infty$ in its diagonal
embedding, projects to a masa of $\Cc(l^2)$ which is isomorphic to
$l^\infty/c_0$. Call this the {\it standard masa} of $\Cc(l^2)$.

Assuming CH, Huruya \cite{Hur} proved that there exists a bounded linear
map from some $C(X)$ into $l^\infty/c_0 \cong C(\beta \omega - \omega)$
which cannot be expressed as a linear combination of positive linear
maps; moreover, this is still true if $l^\infty/c_0$ is identified with
the standard masa of $\Cc(l^2)$ and the positive linear maps are allowed
to go from $C(X)$ to $\Cc(l^2)$. It seems likely that other independence
phenomena involving Banach space properties of $l^\infty/c_0$ (e.g.,
see \cite{DR}) should have analogs for $\Cc(l^2)$.

Anderson \cite{A2} showed that CH implies there are other masas of
$\Cc(l^2)$ besides those of the form $\pi(\Ac)$ for $\Ac$ a masa of
$\Bc(l^2)$. He did this by proving that there is a masa in $\Cc(l^2)$
which contains an uncountable family of commuting projections $\{p_\alpha\}$
such that every element of $\Cc(l^2)$ commutes with all but countably
many $p_\alpha$. He also proved that assuming CH, for every countable
set of states on $\Cc(l^2)$ there is a masa $\Ac$ such that the
restriction of each state to $\Ac$ is pure.

Anderson \cite{A1} also used CH to construct a masa $\Ac$ of $\Cc(l^2)$
with the property that every pure state on $\Ac$ has only one (necessarily
pure) state extension to $\Cc(l^2)$; it is unclear whether this
masa lifts to a masa of $\Bc(l^2)$. This relates to the problem of
classifying the pure states on $\Cc(l^2)$, which is a major unsolved
problem in C*-algebra theory.

The centrality of this pure state problem arises from the fact that
every pure state on $\Cc(l^2)$ is also, after composition with $\pi$, a
pure state on $\Bc(l^2)$. As any pure state on $\Bc(l^2)$ is either a
vector state or else it factors through $\Cc(l^2)$, this means that the
pure states on $\Cc(l^2)$ correspond precisely to the nontrivial pure
states on $\Bc(l^2)$ (much as the points of $\beta\omega - \omega$
correspond precisely to the nontrivial ultrafilters over $\omega$).

Let $\Uc$ be a free ultrafilter over $\omega$ and let $(e_n)$ be
the standard orthonormal basis of $l^2$. Then the map
$A \mapsto \lim_\Uc \langle Ae_n, e_n\rangle$ is a pure state on
$\Bc(l^2)$ \cite{A3} which annihilates the compact operators and
therefore also defines a pure state on $\Cc(l^2)$. Anderson \cite{A5}
conjectured that every nontrivial pure state on $\Bc(l^2)$ is of
this form with respect to some orthonormal basis of $l^2$.

A more famous question called the {\it Kadison-Singer problem}
\cite{KS} asks whether for every free ultrafilter $\Uc$
the pure state $a \mapsto \lim_\Uc a_n$
on $l^\infty$ (i.e., evaluation of the sequence $a = (a_n) \in
l^\infty \cong C(\beta\omega)$ at $\Uc \in \beta\omega$) has
only one extension to a pure state of $\Bc(l^2)$ when $l^\infty$
is diagonally embedded in $\Bc(l^2)$. This problem has attracted
a great deal of attention because it has important implications
in a variety of disparate areas \cite{CT}. It seems unlikely to be
independent of ZFC because it can be restated in simple combinatorial
terms involving partitioning finite sets of vectors in $\Cb^n$ \cite{W3}
(and this formulation is arithmetical). Anderson's conjecture mentioned
above is a more likely candidate
for independence. However, partial results on the Kadison-Singer
problem have used CH: Reid \cite{Rei} showed that every rare
ultrafilter over $\omega$ has the unique extension property
(and noted that CH implies the existence of rare ultrafilters),
and Anderson \cite{A4} directly used CH to construct an ultrafilter
that is not rare but still has the unique extension property.

If projections in $\Cc(l^2)$ are analogous to clopen subsets of
$\beta\omega - \omega$, then pure states on $\Cc(l^2)$ (or
equivalently, nontrivial pure states on $\Bc(l^2)$) may be
seen as analogous to free ultrafilters over $\omega$. There
should be many interesting issues to pursue here, besides
those mentioned above.



\bigskip
\bigskip
\begin{thebibliography}{aaaaaaaa}

\bibitem [1]{AAP1}
C.\ A.\ Akemann, J.\ Anderson, and G.\ K.\ Pedersen,
Excising states of C*-algebras, {\it Can.\ J.\ Math.\
\bf 38} (1986), 1239-1260.

\bibitem [2]{AAP2}
{---------},
Diffuse sequences and perfect C*-algebras, {\it Trans.\ Amer.\ Math.\ Soc.
\bf298} (1986), 747-762.

\bibitem [3]{AP}
C.\ A.\ Akemann and G.\ K.\ Pedersen, Ideal perturbations of elements
in C*-algebras, {\it Math.\ Scand.\ \bf 41} (1977), 117-139.

\bibitem [4]{AS}
C.\ A.\ Akemann and F.\ W.\ Shultz, {\it Perfect C*-algebras},
Mem.\ Amer.\ Math.\ Soc.\ No.\ 326, 1985.

\bibitem [5]{AW}
C.\ Akemann and N.\ Weaver, Consistency of a counterexample to Naimark's
problem, {\it Proc.\ Natl.\ Acad.\ Sci.\ USA \bf 101} (2004), 7522-7525.

\bibitem [6]{AGPS}
S.\ Albeverio, D.\ Guido, A.\ Ponosov, and S.\ Scarlatti, Singular
traces and compact operators, {\it J.\ Funct.\ Anal.\ \bf 137}
(1996), 281-302.

\bibitem [7]{A1}
J.\ Anderson, A maximal abelian subalgebra of the Calkin algebra
with the extension property, {\it Math.\ Scand.\ \bf 42} (1978),
101-110.

\bibitem [8]{A2}
{---------}, Pathology in the Calkin algebra, {\it J.\ Operator
Theory \bf 2} (1979), 159-167.

\bibitem [9]{A3}
{---------}, Extreme points in sets of positive linear maps on
$\Bc(H)$, {\it J.\ Funct.\ Anal.\ \bf 31} (1979), 195-217.

\bibitem [10]{A4}
{---------}, Extensions, restrictions, and representations of
states on C*-algebras, {\it Trans.\ Amer.\ Math.\ Soc.\ \bf 249}
(1979), 303-329.

\bibitem [11]{A5}
{---------}, A conjecture concerning the pure states of $\Bc(H)$ and
a related theorem, in {\it Topics in Modern Operator Theory}, pp. 27-43,
Birkha\"user, 1981. 

\bibitem [12]{Ara}
H.\ Araki, {\it Mathematical Theory of Quantum Fields}, Oxford University
Press, 1999.

\bibitem [13]{Arc}
R.\ Archbold, On perfect C*-algebras, {\it Proc.\ Amer.\ Math.\ Soc.\
\bf 97} (1986), 413-417.

\bibitem [14]{Blk}
B.\ E.\ Blackadar, Weak expectations and nuclear C*-algebras,
{\it Indiana Univ.\ Math.\ J.\ \bf 27} (1978), 1021-1026.

\bibitem [15]{Bla}
A.\ Blass, Near coherence of filters II: applications to operator
ideals, the Stone-\v{C}ech remainder of a half-line, order ideals
of sequences, and slenderness of groups, {\it Trans.\ Amer.\ Math.\
Soc.\ \bf 300} (1987), 557-581.

\bibitem [16]{BS}
A.\ Blass and S.\ Shelah, Near coherence of filters III: a simplified
consistency proof, {\it Notre Dame J.\ Formal Logic \bf 30} (1989),
530-538.

\bibitem [17]{BW}
A.\ Blass and G.\ Weiss, A characterization and sum decomposition
for operator ideals, {\it Trans.\ Amer.\ Math.\ Soc.\ \bf 246}
(1978), 407-417.

\bibitem [18]{BR}
O.\ Bratteli and D.\ W.\ Robinson, {\it Operator Algebras and
Quantum Statistical Methanics 2: Equilibrium States, Models in
Quantum Statistical Mechanics (Second edition)}, Springer-Verlag, 1997.

\bibitem [19]{CT}
P.\ G.\ Casazza and J.\ C.\ Tremain, The Kadison Singer problem in
mathematics and engineering, {\it Proc.\ Nat.\ Acad.\ Sci.\ USA \bf 103}
(2006), 2032-2039.

\bibitem [20]{Cnn}
A.\ Connes, A survey of foliations and operator algebras, in
{\it Operator Algebras and Applications, Part I}, pp.\ 521-628,
American Mathematical Society, 1982.

\bibitem [21]{Con}
J.\ B.\ Conway, {\it A Course in Functional Analysis} (second edition),
Springer-Verlag, 1990.

\bibitem [22]{DW}
H.\ G.\ Dales and W.\ H.\ Woodin, {\it An Introduction to Independence
for Analysts}, Cambridge University Press, 1987.

\bibitem [23]{Dav}
K.\ R.\ Davidson, Similarity and compact perturbations of nest algebras,
{\it J.\ Reine Angew.\ Math.\ \bf 348} (1984), 72-87.

\bibitem [24]{Dix}
J.\ Dixmier, Anneaux d'op\'erateurs et repr\'esentations des groupes,
{\it S\'eminaire Bourbaki, Vol. 1}, Exp.\ No.\ 40, pp.\ 331-336,
Soc.\ Math.\ France, 1995.

\bibitem [25]{Dow}
A.\ Dow, On ultrapowers of Boolean algebras, {\it Topology Proc.\ \bf 9}
(1984), 269-291.

\bibitem [26]{DR}
L.\ Drewnowski and J.\ W.\ Roberts, On the primariness of the Banach
space $l^\infty/C_0$, {\it Proc.\ Amer.\ Math.\ Soc.\ \bf 112}
(1991), 949-957.

\bibitem [27]{GH}
L.\ Ge and D.\ Hadwin, Ultraproducts of C*-algebras, in
{\it Operator Theory: Advances and Applications}, pp.\ 305-326,
Birkh\"auser, 2001.

\bibitem [28]{Gli}
J.\ Glimm, Type I C*-algebras, {\it Ann.\ of Math.\ \bf 73} (1961),
572-612.

\bibitem [29]{Had}
D.\ Hadwin, Maximal nests in the Calkin algebra, {\it Proc.\ Amer.\
Math.\ Soc.\ \bf 126} (1998), 1109-1113.

\bibitem [30]{HKL}
L.\ A.\ Harrington, A.\ S.\ Kechris, and A.\ Louveau,
A Glimm-Effros dichotomy for Borel equivalence relations,
{\it J.\ Amer.\ Math.\ Soc.\ \bf 3} (1990), 903-928.

\bibitem [31]{HO}
T.\ Hinokuma and M.\ Ozawa, Conversion from nonstandard matrix algebras
to standard factors of type $II_1$, {\it Ill.\ J.\ Math.\ \bf 37}
(1993), 1-13.

\bibitem [32]{Hur}
T.\ Huruya, Decompositions of linear maps into nonseparable C*-algebras,
{\it Publ.\ Res.\ Inst.\ Math.\ Sci.\ \bf 21} (1985), 645-655.

\bibitem [33]{Jec}
T.\ Jech, {\it Set Theory}, Springer-Verlag, 2003.

\bibitem [34]{JP}
B.\ E.\ Johnson and S.\ K.\ Parrott, Operators commuting with a von
Neumann algebra modulo the set of compact operators, {\it J.\ Funct.\
Anal.\ \bf 11} (1972), 39-61.

\bibitem [35]{Jon}
V.\ F.\ R.\ Jones, {\it Subfactors and Knots}, American Mathematical
Society, 1991.

\bibitem [36]{Jor}
P.\ E.\ T.\ Jorgensen, {\it Analysis and Probability: Wavelets, Signals,
Fractals}, Springer-Verlag, 2006.

\bibitem [37]{KR}
R.\ V.\ Kadison and J.\ R.\ Ringrose, {\it Fundamentals of the Theory of
Operator Algebras, Vol. I}, AMS, 1997.

\bibitem [38]{KS}
R.\ V.\ Kadison and I.\ M.\ Singer, Extensions of pure states,
{\it Amer.\ J.\ Math.\ \bf 81} (1959), 383-400.

\bibitem [39]{Kas}
G.\ G.\ Kasparov, Equivariant $KK$-theory and the Novikov conjecture,
{\it Invent.\ Math.\ \bf 91} (1988), 147-201.

\bibitem [40]{Kun}
K.\ Kunen, {\it Set Theory: An Introduction to Independence Proofs},
North-Holland, 1980.

\bibitem [41]{Lau}
K.\ B.\ Laursen, Continuity of homomorphisms from C*-algebras into
commutative Banach algebras, {\it J.\ London Math.\ Soc.\ \bf 36}
(1987), 165-175.

\bibitem [42]{Mau}
R.\ D.\ Mauldin, A representation theorem for the second dual of $C[0,1]$,
{\it Studia Math.\ \bf 46} (1973), 197-200.

\bibitem [43]{Mun}
D.\ Mundici, Simple Bratteli diagrams with a G\"odel-incomplete
C*-equivalence problem, {\it Trans.\ Amer.\ Math.\ Soc.\ \bf 356}
(2004), 1937-1955.

\bibitem [44]{N1}
M.\ A.\ Naimark, Rings with involutions (Russian), {\it Uspehi Matem.\
Nauk (N.S.) \bf 3} (1948), 52-145.

\bibitem [45]{N2}
{---------}, On a problem of the theory of rings with involution (Russian),
{\it Uspehi Matem.\ Nauk (N.S.) \bf 6} (1951), 160-164.

\bibitem [46]{Ped}
G.\ K.\ Pedersen, {\it C*-Algebras and Their Automorphism Groups},
Academic Press, 1979.

\bibitem [47]{Rei}
G.\ A.\ Reid, On the Calkin representations, {\it Proc.\ London
Math.\ Soc.\ \bf 23} (1971), 547-564.

\bibitem [48]{Ros}
A.\ Rosenberg, The number of irreducible representations of simple
rings with no minimal ideals, {\it Amer.\ J.\ Math.\ \bf 75} (1953),
523-530.

\bibitem [49]{Run}
V.\ Runde, The structure of discontinuous homomorphisms from
non-commutative C*-algebras, {\it Glasgow Math.\ J.\ \bf 36}
(1994), 209-218.

\bibitem [50]{Sak}
S.\ Sakai, {\it C*-Algebras and W*-Algebras}, Springer-Verlag, 1971.

\bibitem [51]{S1}
S.\ Shelah, {\it Proper Forcing}, Springer-Verlag, 1982.

\bibitem [52]{S2}
{---------}, {\it Classification Theory and the Number of Nonisomorphic
Models}, North-Holland, 1990.

\bibitem [53]{SS}
S.\ Shelah and J.\ Stepr\=ans, PFA implies all automorphisms are trivial,
{\it Proc. Amer. Math. Soc. \bf 104} (1988), 1220-1225.

\bibitem [54]{Tak}
M.\ Takesaki, {\it Theory of Operator Algebras I}, Springer-Verlag,
1979.

\bibitem [55]{vM}
J.\ van Mill, An introduction to $\beta\omega$, in {\it Handbook
of Set-Theoretic Topology}, pp.\ 503-567, North-Holland, 1984.

\bibitem [56]{Vel}
B.\ Veli\v{c}kovi\'c, OCA and automorphisms of $P(\omega)/{\rm fin}$,
{\it Topology Appl.\ \bf 49} (1993), 1-13.

\bibitem [57]{Ves}
J.\ Vesterstr{\o}m, Quotients of finite W*-algebras, {\it J.\
Funct.\ Anal.\ \bf 9} (1972), 322-335.

\bibitem [58]{Wal}
R.\ M.\ Wald, {\it Quantum Field Theory in Curved Spacetime and
Black Hole Thermodynamics}, University of Chicago Press, 1994.

\bibitem [59]{W1}
N.\ Weaver, {\it Mathematical Quantization}, CRC Press, 2001.

\bibitem [60]{W2}
{---------}, A prime C*-algebra that is not primitive, {\it  J.\ Funct.\
Anal.\ \bf 203} (2003), 356-361.

\bibitem [61]{W3}
{---------}, The Kadison-Singer problem in discrepancy theory,
{\it Discrete Math.\ \bf 278} (2004), 227-239.

\bibitem [62]{Wof}
E.\ Wofsey, Set theory and projections in the Calkin algebra, manuscript.

\end{thebibliography}
\end{document}